\documentstyle[psfig,amssymb]{llncs}
\begin{document}

\newcommand{\beq}{\begin{equation}}
\newcommand{\eeq}{\end{equation}}
\newcommand{\bea}{\begin{eqnarray*}}
\newcommand{\eea}{\end{eqnarray*}}
\newcommand{\eps}{\epsilon}
\newcommand{\Z}{{\Bbb Z}}
\newcommand{\back}{\backslash}
\newcommand{\cho}{\left(\!\begin{array}{c}}
\newcommand{\ose}{\end{array}\!\right)}
\newcommand{\mat}{\left(\!\begin{array}{cc}}
\newcommand{\rix}{\end{array}\!\right)}
\newcommand{\ising}{{\rm Ising}}
\newcommand{\domega}{{\rm d}\omega}

\title{Some Polyomino Tilings of the Plane}
\author{Cristopher Moore}
\institute{Santa Fe Institute, 1399 Hyde Park Road,
Santa Fe, New Mexico 87501 {\tt moore@santafe.edu}}
\maketitle

\begin{abstract}  
We calculate the generating functions for the number of tilings of
rectangles of various widths by the right tromino, the $L$ tetromino,
and the $T$ tetromino.  This allows us to place lower bounds on the
entropy of tilings of the plane by each of these.  For the $T$
tetromino, we also derive a lower bound from the solution of the Ising
model in two dimensions.
\end{abstract}

\section{Introduction}

Tilings of the plane are of interest both to statistical physicists
and recreational mathematicians.  While the number of ways to tile a
lattice with dominoes can be calculated by expressing it as a
determinant \cite{kasteleyn,propp}, telling whether a finite collection 
of shapes can tile the plane at all is undecidable
\cite{berger,robinson}.

In the 1994 edition of his wonderful book {\em Polyominoes}
\cite{golomb}, Solomon W. Golomb states that the problem of
determining how many ways a $4 \times n$ rectangle can be tiled by
right trominoes ``appears to be a challenging problem with reasonable
hope of an attainable solution.''  Here we solve this problem, and
several more.  Specifically, we find the generating functions for
tilings of $4 \times n$ and $5 \times n$ rectangles by the right
tromino, and tilings of $4 \times n$ rectangles by the $L$ and $T$
tetrominoes. 

The number of ways of tiling a $m \times n$ rectangle with any finite
collection of shapes, where $m$ is fixed, can be found by calculating
the $n$th power of a transfer matrix whose rows and columns correspond
to the various interface shapes a partial tiling can have.  As $n$
grows large, the number of tilings $N(n)$ grows as $\lambda^n$ where
$\lambda$ is the largest eigenvalue of this matrix.  In general, the
number of interface shapes and therefore the size of the transfer
matrix grows exponentially in $m$.  However, if a small number of
interface shapes suffice to generate all possible tilings, then we can
use a small transfer matrix, and $\lambda$ will then be an algebraic
number of low degree.  In fact, we find that for the problems
considered here the largest matrix we need is $4 \times 4$.

By calculating $\lambda$ for these rectangles, we can place a lower bound
on the entropy per site of tilings of the plane by these polyominoes.
In the case of the $T$ tetromino, we improve this bound by using the
exact solution of the Ising model in two dimensions.

\section{The right tromino}

\begin{figure}
\centerline{\psfig{file=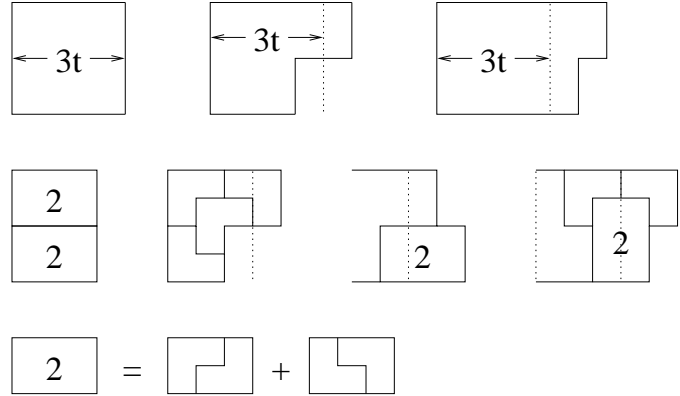,width=3.5in}}
\caption{The three kinds of interfaces for the right tromino on
rectangles of width 4, straight, deep jog, and shallow jog, and the
transitions between them.}
\label{right}
\end{figure}

Consider tilings of a $4 \times n$ rectangle with right trominoes,
where $n$ is a multiple of 3.  The number of interfaces seems
potentially large.  However, it turns out that we only need to think
about three, so that a $3 \times 3$ matrix will suffice.  These are
shown in Figure~\ref{right}, with the various kinds of transitions
that can take place as the tiling grows to the right.

For these three interfaces, which we call `straight,' `deep jog,' and
`shallow jog,' we define $N(t)$, $N_1(t)$ and $N_2(t)$ respectively as
the number of tilings when there are $n=3t$ columns to the left of the
dotted line.  Obviously $N_1(t)$ and $N_2(t)$ stay the same if we
count tilings of their vertical reflections instead.  Our goal is to
find $N(t)$.  To express these, we write them as generating functions,
\[ G(z) = \sum_t N(t) \,z^t \]
and similarly for $G_1(z)$ and $G_2(z)$.  

From Figure~\ref{right}, we see that a straight interface can make a
transition to itself in four ways, each of which increases $t$ by 1.
It has one transition each to a deep jog and its reflection, each of
which has one transition back to a straight, and both of these
increase $t$ by 1.  Deep and shallow jogs of a given orientation have
two transitions in each direction; as we have defined them a
shallow-to-deep transition increments $t$ but a deep-to-shallow
transition does not.  A shallow jog has two transitions to its
reflection, incrementing $t$.  Finally, $N(0)=1$ since there is
exactly one way to tile a $4 \times 0$ rectangle.

This gives us a $3 \times 3$ transfer matrix, and the system of linear
equations
\bea G & = & 1 + 4 z G + 2 z G_1 \\
   G_1 & = & z G + 2 z G_2 \\
   G_2 & = & 2 G_1 + 2 z G_2 \eea
The solution for $G$ is
\[ G(z) = \frac{1 - 6 z}{1 - 10 z + 22 z^2 + 4 z^3} \]
The coefficients $N(t)$ of the first few terms of $G$'s Taylor expansion are
\[ 1, 4, 18, 88, 468, 2672, 16072, 100064, 636368, 4097984, 26579488, \ldots
\]
giving the number of ways to tile rectangles of size $4 \times 0$, $4
\times 3$, $4 \times 6$, and so on.

\begin{figure}
\centerline{\psfig{file=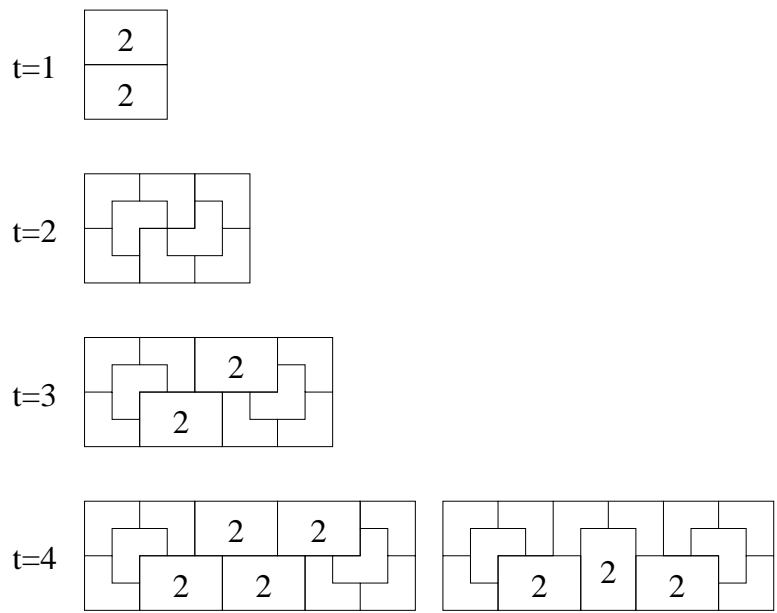,width=4in}}
\caption{Fault-free rectangles of right trominoes of width 4 and
length 3, 6, 9, and 12.}
\label{rightcount}
\end{figure}

To gain a better understanding of this series, we concentrate on {\em
fault-free} rectangles, those whose only straight interfaces are those
at the left and right end.  The first few of these are shown in
Figure~\ref{rightcount}.  These always begin and end with transitions
from a straight edge to a deep jog and back, except for $t=1$ where we
have a pair of $2 \times 3$ rectangles.  By disallowing transitions
back to a straight interface except at the end, we find that the
linear equations for fault-free rectangles are
\bea G' & = & 4 z + 2 z G'_1 \\
   G'_1 & = & z + 2 z G'_2 \\
   G'_2 & = & 2 G'_1 + 2 z G'_2 \eea
This has the solution
\[ G'(z) = 2 z \,\frac{2 - 11 z - 2 z^2}{1 - 6 z} \]
The Taylor series of $G'$ tells us that the number $N'$ of fault-free
tilings is 
\[ 4, 2, 8, 48, 288, 1728, \ldots \]
which we can write as
\[ N'(t) = \left\{ \begin{array}{ll}
   4 & \quad \mbox{if } t=1 \\
   2 & \quad \mbox{if } t=2 \\
   8 \cdot 6^{t-3} & \quad \mbox{if } t \ge 3 
\end{array} \right. \]
Since any tiling consists of a concatenation of fault-free tilings,
we have
\[ G(z) = 1 + G'(z) + G'(z)^2 + G'(z)^3 + \cdots = \frac{1}{1-G'(z)} \]
which the reader can verify.  Note that we have to define
$G'(0)=N'(0)=0$ in order for this formula to work.

The asymptotic growth of $N(t)$ is the reciprocal of the radius of
convergence of $G$'s Taylor series.  Thus $N(t) \propto \lambda^t$
where $\lambda$ is the largest positive root of 
\[ \lambda^3 - 10 \lambda^2 + 22 \lambda + 4 = 0 \]
This gives
\[ \lambda = \frac{2}{3} \, 
   \left( 5 + \sqrt{34}\,\cos \frac{\theta}{3} \right) \]
where $\pi/2 < \theta < \pi$ and
\[ \tan \theta = - \frac{3}{11} \, \sqrt{\frac{519}{2}} \]
Numerically, we have
\[ \lambda = 6.54560770847481152029 \cdots \]
It would be nice to find a simpler expression for $\lambda$, perhaps
using the decomposition into fault-free rectangles.

\begin{figure}
\centerline{\psfig{file=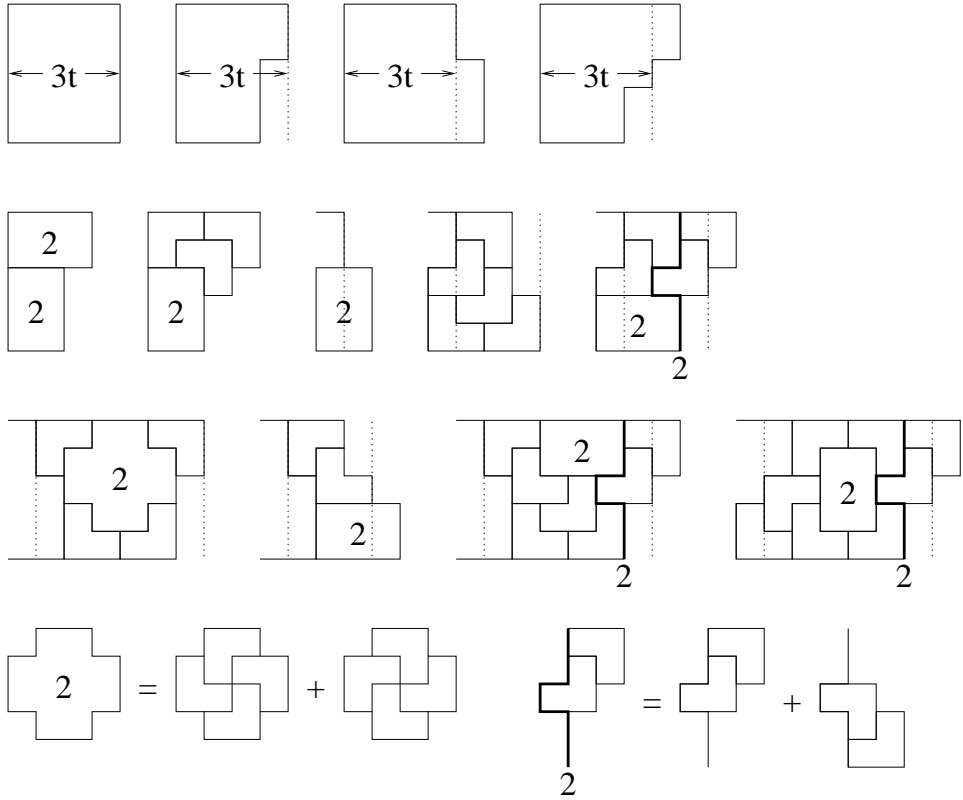,width=5in}}
\caption{The four kinds of interfaces for the right tromino on
rectangles of width 5, straight, big jog, little jog, and slope,
and the transitions between them.}
\label{right5}
\end{figure}

For $5 \times 3t$ rectangles, four interface shapes suffice.  We call
these straight, little jog, big jog, and slope, and their transitions
are shown in Figure~\ref{right5} (reflections and reversals are not
shown).  If their generating functions are $G$, $G_1$, $G_2$ and $G_3$
respectively, our $4 \times 4$ transfer matrix is
\bea G & = & 1 + 8 z G_2 + 4 z G_3 \\
   G_1 & = & 4 z G + z G_1 + 2 z^2 G_2 + (2 z + 4 z^2) G_3 \\
   G_2 & = & 2 G_1 + z G_2 + 4 z G_3 \\
   G_3 & = & 2 z G + 4 z G_1 + (2 z + 4 z^2) G_2 + 4 z^2 G_3 \eea
whose solution is
\[ G(z) = \frac{1 - 2 z - 31 z^2 - 40 z^3 - 20 z^4}
               {1 - 2 z - 103 z^2 + 280 z^3 + 380 z^4} \] 
The coefficients $N(t)$ of the first few terms of $G$'s Taylor
expansion are
\[ 1, 0, 72, 384, 8544, 76800, 1168512, 12785664, 
   170678784, 2014648320, 25633231872, \ldots \] 
Note that there are no $5 \times 3$ tilings.

\begin{figure}
\centerline{\psfig{file=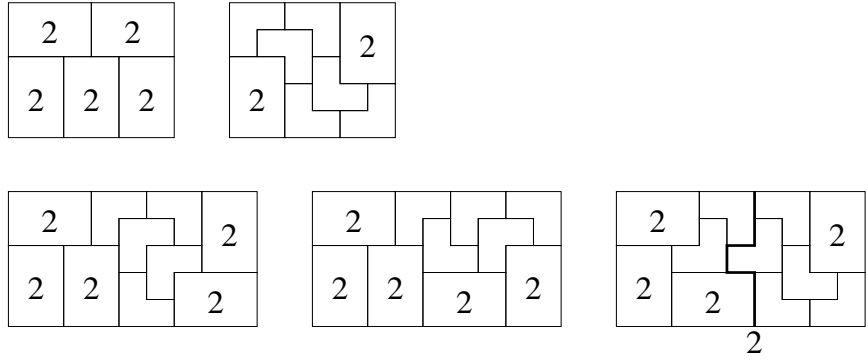,width=4.5in}}
\caption{Fault-free rectangles of right trominoes of width 5 and
length 6 and 9.}
\label{right5count}
\end{figure}

The generating functions for fault-free tilings obey
\bea G' & = & 8 z G'_2 + 4 z G'_3 \\
   G'_1 & = & 4 z + z G'_1 + 2 z^2 G'_2 + (2 z + 4 z^2) G'_3 \\
   G'_2 & = & 2 G'_1 + z G'_2 + 4 z G'_3 \\
   G'_3 & = & 2 z + 4 z G'_1 + (2 z + 4 z^2) G'_2 + 4 z^2 G'_3 \eea
and so
\[ G'(z) = 24 z^2 \,\frac{3 + 10 z + 15 z^2}
                         {1 - 2 z - 31 z^2 - 40 z^3 - 20 z^4} \] 
The first few $N'(t)$, starting with $t=2$, are then 
\[ 72, 384, 3360, 21504, 163968,
1136640, 8283648, 58791936, 423121920, \ldots \]
These are shown for $t=2$ and $t=3$ in Figure~\ref{right5count}.

The asymptotic growth of the number of fault-free rectangles is $N'(t)
\propto \lambda'^t$ where $\lambda'$ is the largest positive root of 
\[ \lambda'^4 - 2 \lambda'^3 - 31 \lambda'^2 - 40 \lambda' - 20 = 0 \]
which has a rather complicated closed-form expression which we will
not reproduce here.  Numerically, 
\[ \lambda' = 7.16235536278185348653 \cdots \]
For all rectangles, including faulty ones, we have $N(t) \propto
\lambda^t$ where $\lambda$ is the largest positive root of
\[ \lambda^4 - 2 \lambda^3 - 103 \lambda^2 - 280 \lambda - 380 = 0 \]
This is
\[ \lambda = 
  \frac{1}{2}  
+ {\frac{{\sqrt{627\,y^\frac{1}{3} - 3\,y^\frac{2}{3} - 13107}}}
   {6\,y^\frac{1}{6}}} + 
   {\frac{{\sqrt{4369 + 418\,{y^\frac{1}{3}} + {y^\frac{2}{3}} + 
          2304\,{\sqrt{\frac{3\,y}
                 {209\,y^\frac{1}{3} - y^\frac{2}{3} - 4369}}}}}}
       {2\,\sqrt{3}\,y^\frac{1}{6}}}
\]
where
\[ y = 1204327 - 72\,\sqrt{263697405} \]
Numerically, we have
\[ \lambda = 12.36366722455963019234 \cdots \]

\section{$L$ tetrominoes}

\begin{figure}
\centerline{\psfig{file=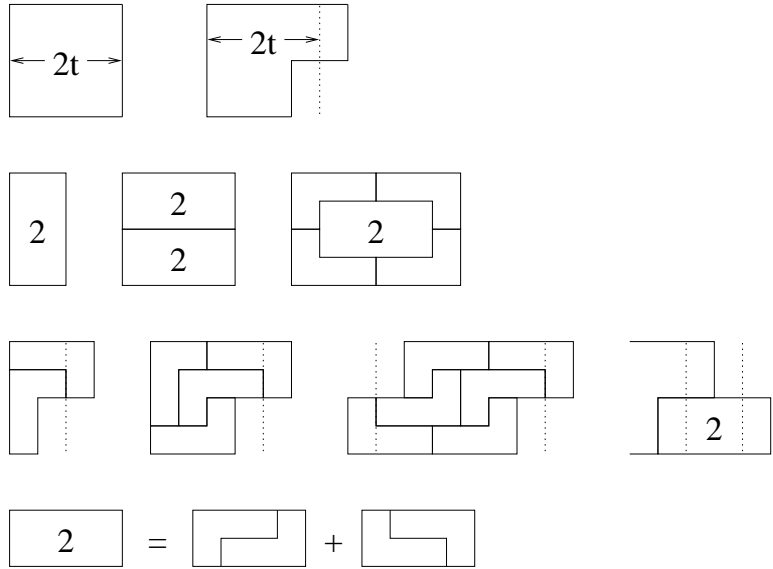,width=4in}}
\caption{The two kinds of interfaces for the $L$ tetromino on
rectangles of width 4, straight and jogged, and the transitions
between them.}
\label{ltet}
\end{figure}

We can perform a similar analysis for tilings of $4 \times n$
rectangles with $L$ tetrominoes.  In fact, here the analysis is
easier.  As Figure~\ref{ltet} shows, we have just two kinds of
interfaces, straight and jogged.  We define $N(t)$ and $N_1(t)$ as the
number of tilings for a straight or jogged interface respectively when
there are $2t$ columns to the left of the dotted line.  Then a
straight has two transitions to itself that increase $t$ by 1, four to
itself that increase $t$ by 2, and two to itself that increase $t$ by
3.  It has two transitions each to a jog and its reflection, one of
which increases $t$ by 1 and the other by 2.  Finally, a jog has two
transitions to another (reflected) jog that increase $t$ by 1, one
transition to itself that increases $t$ by 3, and the two inverse
transitions back to a straight.

Thus we have a $2 \times 2$ transfer matrix, and defining generating
functions as before gives us the linear equations
\bea G & = & 1 + (2 z + 4 z^2 + 2 z^3) G + (2 z + 2 z^2) G_1 \\
   G_1 & = & (z + z^2) G + (2 z + z^3) G_1 \eea
so
\[ G(z) = \frac{1 - 2 z - z^3}
     {1 - 4 z - 2 z^2 + z^3 + 4 z^4 + 4 z^5 + 2 z^6} \]
The coefficients $N(t)$ of the first few terms of $G$'s Taylor
expansion are
\[ 1, 2, 10, 42, 182, 790, 3432, 14914, 64814, 281680, 1224182, \ldots \]
giving the number of ways to tile a rectangle of size $4 \times 0$, $4
\times 2$, $4 \times 4$, etc.

\begin{figure}
\centerline{\psfig{file=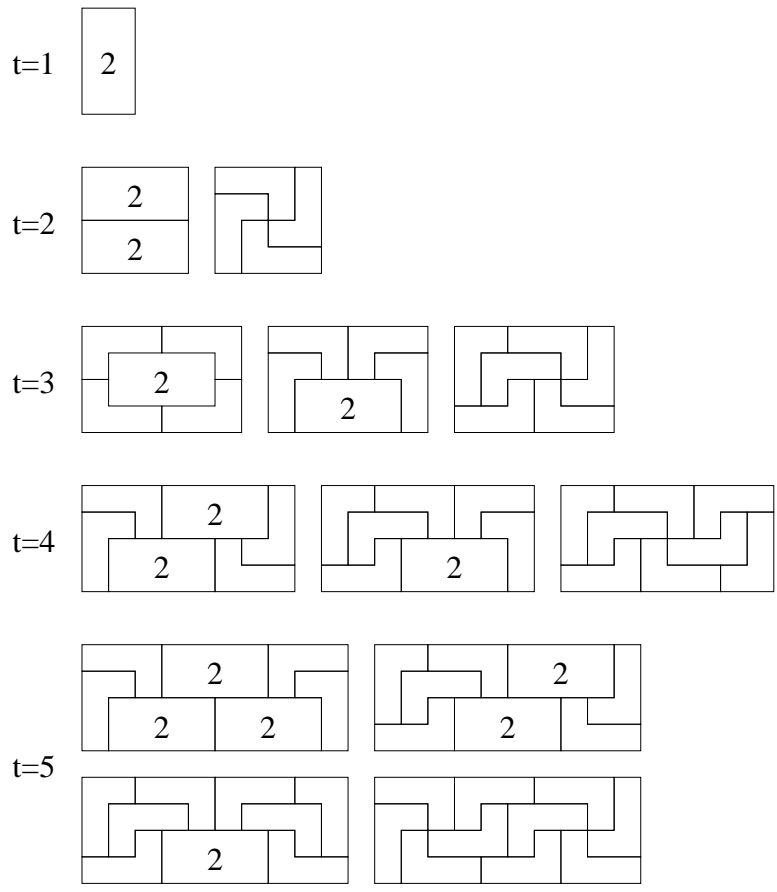,width=4in}}
\caption{Fault-free rectangles of $L$ tetrominoes of width 4 and
length 2, 4, 6, 8, and 10.}
\label{ltetcount}
\end{figure}

As before, we can focus our attention on fault-free rectangles.  For
$t > 3$, any such rectangle consists of transitions from a straight
edge to a jog, and transitions between jogs in between.  For $t=1, 2,
3$ we have an additional 2, 4, and 2 fault-free rectangles
respectively as shown in Figure~\ref{ltet}.  The linear equations for
the fault-free generating functions are thus
\bea G' & = & 2 z + 4 z^2 + 2 z^3 + 2 (z + z^2) G'_1 \\
   G'_1 & = & z + z^2 + (2 z + z^3) G'_1 \eea
This gives
\[ G'(z) = \frac{2 z\,(1 + z)^2\,(1 - z - z^3)}{1 - 2 z - z^3} \]
whose Taylor expansion tells us that the number $N'$ of fault-free
tilings is
\[ 2, 6, 10, 18, 38, 84, 186, 410, 904, 1994 \ldots \]
The first few of these are shown in Figure~\ref{ltetcount}.  The
reader can verify that $G(z) = 1/(1-G'(z))$.

The asymptotic growth of the number of fault-free rectangles is $N'(t)
\propto \lambda'^t$ where $\lambda'$ is the largest root of
\[ \lambda'^3 - 2 \lambda'^2 - 1 = 0 \]
which is
\[ \lambda' = \frac{1}{3} \left(
  2 + \sqrt[3]{\frac{43}{2} - \frac{3\,\sqrt{177}}{2}}
    + \sqrt[3]{\frac{43}{2} + \frac{3\,\sqrt{177}}{2}} \right) \]
or numerically
\[ \lambda' = 2.20556943040059031170 \cdots \]
For all rectangles, including faulty ones, we have $N(t) \propto
\lambda^t$ where $\lambda$ is the largest root of 
\[ \lambda^6 - 4 \lambda^5 - 2 \lambda^4 + \lambda^3 + 4 \lambda^2 
   + 4 \lambda + 2 = 0 \]
This appears not to have a closed form, but it is approximately
\[ \lambda = 4.34601641141649282849 \cdots \]

\section{$T$ tetrominoes}

\begin{figure}
\centerline{\psfig{file=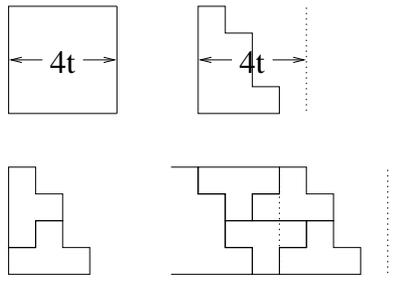,width=2in}}
\caption{The two kinds of interfaces for the $T$-tetromino, straight
and jagged, and the transitions between them.}
\label{ttet}
\end{figure}

\begin{figure}
\centerline{\psfig{file=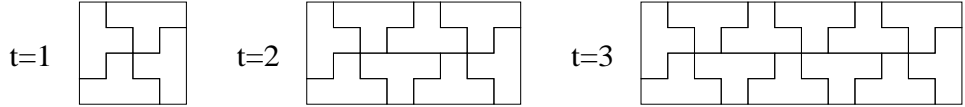,width=5in}}
\caption{The $T$ tetromino has two fault-free rectangles of each size.}
\label{ttetcount}
\end{figure}

Of these three polyominoes, finding the number of tilings of a $4
\times n$ rectangle is easiest for the $T$ tetromino.  In fact, such
tilings only exist if $n$ is a multiple of 4.  Figure~\ref{ttet} shows
two kinds of interfaces, straight and jagged, with $4t$ columns to the
left of the dotted line.  A straight edge must make a transition
either to a jag or its reflection, incrementing $t$.  A jag can make a
transition back to straight, leaving $t$ fixed, or to itself,
incrementing $t$.  Thus the generating functions obey 
\bea G & = & 1 + 2 G_1 \\
   G_1 & = & z G + z G_1 \eea
This gives
\[ G(z) = \frac{1 - z}{1 - 3 z} \]
whose Taylor expansion tells us that
\[ N(t) = \left\{ \begin{array}{ll} 
   1 & \quad \mbox{if } t = 0 \\
   2 \cdot 3^{t-1} & \quad \mbox{if } t > 0
\end{array} \right. \]
Thus we can find a closed form for $N(t)$, unlike in the previous two
cases.  There are exactly two fault-free rectangles for each $t$, as
shown in Figure~\ref{ttetcount}.  Thus the number of tilings is $2
\cdot 3^{t-1}$ since there are two choices of initial fault-free
rectangle and three choices for each increment of $t$, namely either 
continuing the current fault-free rectangle, or ending it and starting
one of two new ones.

\begin{figure}
\centerline{\psfig{file=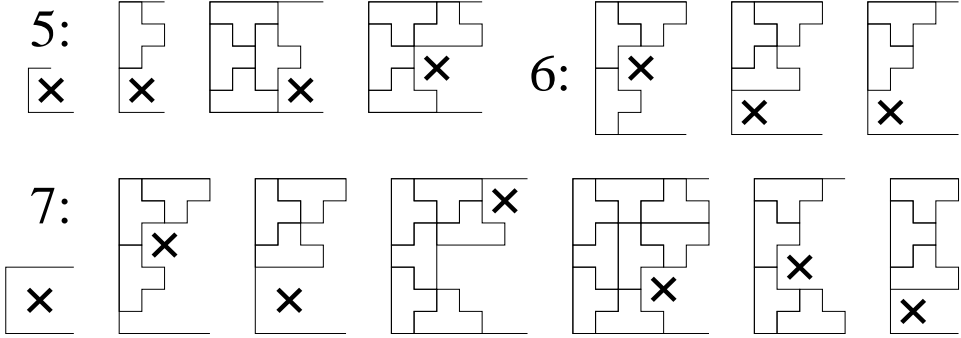,width=5in}}
\caption{Proofs that the $T$ tetromino cannot tile rectangles of width
5, 6, or 7.  Attempted tilings always lead to one of the un-tileable notches
shown on the left, or to an even more obviously un-tileable space.}
\label{ttetcant}
\end{figure}

It is easy to show that the $T$ tetromino cannot tile any rectangles
of width 5, 6, or 7.  Proofs are shown in Figure~\ref{ttetcant}.
More generally, Walkup \cite{walkup} showed that a rectangle can be
tiled with $T$ tetrominoes if and only if its length and width are
both multiples of 4.  

\section{Tilings of the plane}

Since the plane can be divided into rows or columns, the entropy per site 
of tilings of rectangles of a given width serves as a lower bound for the 
entropy per site of tilings of the plane.  Specifically, if the number of 
ways to tile a $m \times n$ rectangle is $N(m,n)$, we define the entropy 
$\sigma$ as
\[ \sigma = \lim_{m,n \to \infty} \frac{\ln N(m,n)}{mn} \]
so that
\[ N(m,n) \propto e^{\sigma mn} \]
for large rectangles.  For the right tromino, our analysis of $5
\times n$ rectangles gives
\[ \sigma \ge \frac{1}{15} \ln \lambda = 0.1676508 \cdots \]
where the factor of 15 comes from the fact that each increment of $t$
adds 15 sites to the rectangle.  Similarly, for the $L$ tetromino on
$4 \times n$ rectangles, we have
\[ \sigma \ge \frac{1}{8} \ln \lambda = 0.183657 \cdots \]
Better lower bounds could be obtained by looking at wider rectangles.
Note that if there are any tilings we failed to see, this only
improves these lower bounds, since additional transitions can only
increase $\lambda$.

\begin{figure}
\centerline{\psfig{file=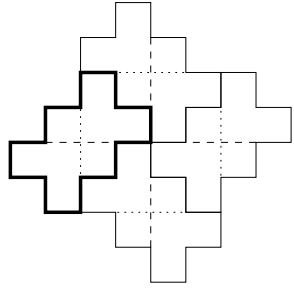,width=1.5in}}
\caption{Since this eight-cell shape can tile the plane and in turn
can be tiled in two ways, the entropy of tilings of the plane with $T$
tetrominoes is at least $\frac{1}{8} \ln 2$.}
\label{ttile}
\end{figure}

For the $T$ tetromino, the lower bound we get from $4 \times n$
rectangles, $\frac{1}{16} \ln 3 \approx 0.06866$, is not as good as
the bound $\frac{1}{8} \ln 2 \approx 0.08664$ which we can derive by
noting that the 8-cell shape in Figure~\ref{ttile} can tile the plane,
and in turn can be tiled in two ways by the $T$.  We can get a better
bound as follows.

\begin{figure}
\centerline{\psfig{file=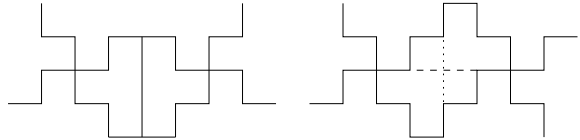,width=3in}}
\caption{The space between two fylfots 4 sites apart can be tiled in 1
way if they are of opposite type, and 2 ways if they are the same.}
\label{fylfot}
\end{figure}

Suppose that the boundaries between tiles at points $(4x, 4y)$ are
clockwise or counterclockwise fylfots as shown in Figure~\ref{fylfot}.
Whenever neighboring fylfots have opposite orientations, there is only
one way to tile the space between them, but if they have the same
orientation the space between them can be tiled in two ways as in
Figure~\ref{ttile}.  (Note that this is also the source of two of the
choices in the $4 \times n$ tilings above.)  If we define the two
orientations as $+1$ and $-1$, then we can sum over all configurations
of fylfots $\{s_i = \pm 1\}$, with each one giving us 1 choice for
pairs of unlike neighbors and 2 choices for like ones.  Thus we have a
lower bound of
\[ N(m,n) = \sum_{\{s_i = \pm 1\}} \prod_{ij} \left\{ \begin{array}{ll}
   1 & \mbox{if } s_i s_j = -1 \\
   2 & \mbox{if } s_i s_j = +1 \end{array} \right\} \]
where the product is over all pairs of nearest neighbors, and the sum
is over the configurations of a lattice of size $\frac{m}{4} \times
\frac{n}{4} = \frac{mn}{16}$.  We can rewrite this as
\bea N(m,n) & = & \sum_{\{s_i = \pm 1\}} \prod_{ij} 
   \,e^{\frac{1}{2}\,(s_i s_j + 1) \ln 2} \\
& = & 2^\frac{mn}{16} \sum_{\{s_i = \pm 1\}} \prod_{ij} 
   \,e^{(\frac{1}{2} \ln 2) \,s_i s_j} \eea
(note that there are $\frac{mn}{8}$ edges in the fylfot lattice).

We now note that the latter sum is the partition function of an
antiferromagnetic Ising model with $\beta = \frac{1}{2} \ln 2$.  We
can transform this to a ferromagnetic model by negating the $s_i$ on
one of the checkerboard sublattices.  Using the exact solution of the
Ising model in two dimensions \cite{mark}, we then have
\[ \sigma \ge \frac{1}{16} \,(\ln 2 + \sigma_\ising) \]
where 
\[ \sigma_\ising = \ln 2 + \frac{1}{2} \frac{1}{(2\pi)^2} 
  \int_0^{2\pi} \!\!\!\int_0^{2\pi} \domega_1 \,\domega_2
    \,\ln \left[
      \cosh^2 2\beta - (\sinh 2\beta)(\cos \omega_1 + \cos \omega_2) \right]
\]
For $\beta = \frac{1}{2} \ln 2$, where $\cosh^2 2\beta = \frac{25}{16}$
and $\sinh 2\beta = \frac{3}{4}$, we have 
\[ \sigma_\ising = 0.8270269567 \cdots \]
and so
\[ \sigma \ge 0.09501088358 \cdots \]

For upper bounds, we can generalize an argument given in \cite{propp}
as follows.  If we construct a tiling of the plane with right
trominoes by scanning from top to bottom and left to right, at each
step the first unoccupied site can be given a tromino with only four
different orientations.  Since we have at most $n/3$ such choices, the
entropy is at most $\frac{1}{3} \ln 4 \approx 0.462$.  Similar
considerations give $\frac{1}{4} \ln 8 \approx 0.520$ and $\frac{1}{4}
\ln 4 \approx 0.347$ for the $L$ and $T$ tetromino respectively.
Obviously the gap between our upper and lower bounds is still quite
large.

For dominoes or `dimers,' the entropy is known to be $G/\pi = 0.29156
\cdots$ where $G$ is Catalan's constant \cite{kasteleyn,propp}.  It is
tempting to think that other polyominoes might have exact solutions;
one related model which does is the covering of the triangular lattice
with triangular trimers \cite{trimers}.  However, while the general
problem of covering an arbitrary graph with dimers can be solved in
polynomial time, covering it with triangles is NP-complete \cite{garey},
so generalized versions of the tromino problem are probably hard.

{\bf Acknowledgements.}  I am grateful to Mark Newman, Lauren Ancel, 
Michael Lachmann and Aaron Meyerowitz for helpful conversations, and 
Molly Rose and Spootie the Cat for warmth and friendship.

{\em Note added.}  Meyerowitz \cite{meyerowitz} has calculated the
generating function for the $T$ tetromino in rectangles of width 8.

\end{document}